\documentclass{article}

\usepackage{amssymb, amsmath, amsthm, array}

\usepackage[left=1in, top=1in, right=1in, bottom=1in]{geometry}

\usepackage{mathrsfs}
\newtheorem{theorem}{Theorem}[section]
\newtheorem{lemma}[theorem]{Lemma}

\newtheorem{corollary}[theorem]{Corollary}
\newtheorem{definition}[theorem]{Definition}

\begin{document}

\title{\bf Left and right zero subsemigroups and an elementary approach to the Rees-Sushkevich Theorem}
\author{Julia Maddox \\
        University of Oklahoma}
\date{}

\maketitle

\abstract{In a semigroup $S$ with idempotent $e$, the set of all elements $a$ of $S$ such that $ae=a$ and $ea=e$ forms the maximum left zero subsemigroup of $S$ containing $e$, and the set of all elements $b$ of $S$ such that $be=e$ and $eb=b$ forms the maximum right zero subsemigroup of $S$ containing $e$.

We use these maximal right and left zero subsemigroups along with the maximal subgroups of a semigroup in the determination of maximal right and left subgroups, rectangular band subsemigroups, completely simple subsemigroups, rectangular subgroups, and completely $0$-simple subsemigroups.

An elementary approach of the Rees-Sushkevich Theorem follows. Then we define rectangular $0$-bands and rectangular $0$-groups as completely $0$-simple analogues of rectangular bands and rectangular groups.
}



\section{Introduction}

A left [right] zero semigroup is defined as $ab=a$ [$ab=b$] for any elements $a$ and $b$ belonging to the semigroup. For any element $e$ of a left [right] zero semigroup, $ee=e$, so that $e$ is an idempotent, and the semigroup is a band. 
Let $S$ be a semigroup with an idempotent $e$.
As proven in Theorem \ref{left zero theorem}, this maximum left [right] zero subsemigroup of $S$ containing $e$, labeled $LZ(e)$ [$RZ(e)$], is the set of all $a \in S$ such that $ae=a$ [$ae=e$] and $ea=e$ [$ea=a$].

Let $S$ be a semigroup. For every idempotent $e \in S$, there exists a maximum subgroup of $S$ containing $e$: the subset of $eSe$, the maximum subsemigroup of $S$ containing $e$ as its identity, of units with respect to $e$. We let $H(e)$ be this maximum subgroup of $S$ containing idempotent $e$. 
Any subgroup of semigroup $S$ contains idempotent $e$ if and only if $e$ is its unique identity, and $H'(e)$ is a subgroup of $S$ containing $e$ if and only if $H'(e)$ is a subgroup of $H(e)$.
For distinct idempotents $e$ and $f$, $H(e) \cap H(f) = \emptyset$. Let $E$ be the set of idempotents in $S$, $\{H(e)|e \in E\}$ is the set of all maximal subgroups of $S$. These are also the $\mathscr{H}$-classes of $S$ containing idempotents.

We can use products of $LZ(e)$, $H(e)$, and $RZ(e)$ when defining the structure of some well known constructs in semigroup theory, which contain idempotents. These can also be extended to comparable definitions in semigroups with $0$.

A right [left] group is a left [right] cancellative and right [left] simple semigroup, and it is isomorphic to a product of a group and a right [left] zero semigroup. A rectangular band is a semigroup $S$ such that $aba=a$ for any $a,b \in S$, and it is isomorphic to a product of a left zero semigroup and a right zero semigroup. A completely simple semigroup is a simple semigroup containing a primitive idempotent. A rectangular group is isomorphic to the direct product of a left zero semigroup, a group, and a right zero semigroup (or a group and a rectangular band). A completely $0$-simple semigroup is a $0$-simple semigroup containing a primitive idempotent, and it is isomorphic to a regular Rees matrix semigroup over a group with zero by the Rees-Sushkevich Theorem. In each of these cases, the semigroup class contains an idempotent (or nonzero idempotent in the completely $0$-simple semigroup case). 

For semigroup $S$, if $e \in S$ is an idempotent, then the trivial group $\{e\}$ is a subsemigroup of $S$, which is also a group, left [right] zero semigroup, right [left] group, rectangular band, completely simple semigroup, and rectangular group. For semigroup $S$ with $0$, if $e \in S$ is a nonzero idempotent, $\{e\}^0$ (the adjunction of $\{e\}$ with $0$) is a completely $0$-simple semigroup. Thus the existence of such subsemigroups coincides with idempotents.

For any semigroup $S$ and for any idempotent $e \in S$, we identify all subsemigroups of $S$ containing $e$ for each of these semigroup classes. In this process, we provide an alternate approach to determining the structure of each of these semigroup classes. 
We use left and right zero subsemigroups and subgroups to identify right [left] subgroups, to identify all rectangular band subsemigroups, to identify all completely simple subsemigroups, to identify all rectangular subgroups, to identify all completely $0$-simple subsemgroups and their decompositions in terms of identity and zero sets utilizing the maximality of such types of subsemigroups under certain conditions.

A rectangular band and a rectangular group are completely simple semigroups. In the last sections, we define completely $0$-simple analogues of these motivated by their structures: a rectangular $0$-band and a rectangular $0$-group.

See \cite{clifford},\cite{howie}, and \cite{grillet} for references on fundamental semigroup theory including maximal subgroups, right [left] zero semigroups, right [left] groups, rectangular bands, completely simple semigroups, rectangular groups, and completely $0$-simple semigroups.

\section{Left [Right] Zero Subsemigroups}

Every right zero or left zero semigroup is a band and contains idempotents.
Let $e$ be an idempotent in semigroup $S$. Define $LZ(e)=\{a \in S|\text{$ae=a$ and $ea=e$}\}$, and define $RZ(e)=\{b \in S|\text{$be=e$ and $eb=b$}\}$.

For idempotent $e$ in semigroup $S$, $L_e=\{a \in S|S^{1}a=S^{1}e\}$ is the $\mathscr{L}$-class of $e$, and for any $a \in L_e$, $ae=a$. Also $R_e=\{a \in S|aS^{1}=eS^{1}\}$ is the $\mathscr{R}$-class of $e$, and for any $a \in R_e$, $ea=a$.
Note that for idempotent $e$ in semigroup $S$, $LZ(e)$ is the set of all idempotents in the $\mathscr{L}$-class of $e$, and $RZ(e)$ is the set of all idempotents in the $\mathscr{R}$-class of $e$.

\begin{theorem} \label{left zero theorem}
Let $S$ be a semigroup, and let $e \in S$ be an idempotent.
$LZ(e)$ is the maximum left zero subsemigroup of $S$ containing $e$.
\end{theorem}

Similarly, $RZ(e)$ is the maximum right zero subsemigroup of $S$ containing $e$.

\begin{proof}
Since $e$ is an idempotent, $ee=e$ and $e \in LZ(e)$, so $LZ(e)$ is nonempty. For any $a,b \in LZ(e)$, $ab = (ae)b = a(eb) = ae = a.$
Thus $LZ(e)$ is a left zero semigroup, and it is the maximum left zero subsemigroup containing $e$.
\end{proof}

We can then define an equivalence for left [right] zero semigroups.

\begin{corollary}
Let $S$ be a semigroup. There exists an idempotent $e \in S$ such that $S=LZ(e)$ [$S=RZ(e)$] if and only if $S$ is a left [right] zero semigroup.
\end{corollary}

\begin{corollary}
Let $S$ be a semigroup, and let $e \in S$ be an idempotent. Then $LZ'(e)$ is a left zero subsemigroup of $S$ containing $e$ if and only if $e \in LZ'(e) \subset LZ(e)$.
\end{corollary}

Similarly, $RZ'(e)$ is a right zero subsemigroup of $S$ containing $e$ if and only if $e \in RZ'(e) \subset RZ(e)$.

\begin{lemma} \label{uniqueness lemma}
Let $S$ be a semigroup with idempotent $e$. Any element in $LZ(e)H(e)RZ(e)$ is a unique product, and if $T=LZ'(e)H'(e)RZ'(e)$ for nonempty subsets $LZ'(e) \subset LZ(e)$, $H'(e) \subset H(e)$, and $RZ'(e) \subset RZ(e)$, then $LZ'(e)$, $H'(e)$, and $RZ'(e)$ are unique.
\end{lemma}

\begin{proof}
Let $e_L,f_L \in LZ(e)$, $a,b \in H(e)$, and $e_R,f_R \in RZ(e)$ such that $e_Lae_R=f_Lbf_R$. Then $e(e_Lae_R)e=e(f_Lbf_R)e$, so that $(ee_L)a(e_Re)=(ef_L)b(f_Re)$. Since $ee_L=ef_L=e_Re=f_Re=e$, $eae=ebe$, which implies $a=b$ because $a,b \in H(e)$. Then $e_Lae_R=f_Laf_R$, and $a^{-1}e(e_Lae_R)=a^{-1}e(f_Laf_R)$, so that $(a^{-1}ee_La)e_R=(a^{-1}ef_La)f_R$. Since $ee_L=ef_L=e$, $(a^{-1}ea)e_R=(a^{-1}ea)f_R$, which implies $ee_R=ef_R$ because $a \in H(e)$. Therefore $e_R=f_R$ because $e_R,f_R \in RZ(e)$. Similarly, $e_L=f_L$.
\end{proof}

This also implies that, for any idempotent $e$ is a semigroup $S$, every product decomposition is unique in $H(e)RZ(e)=\{e\}H(e)RZ(e) \subset LZ(e)H(e)RZ(e)$, and every product decomposition is unique in $LZ(e)RZ(e)=LZ(e)\{e\}RZ(e) \subset LZ(e)H(e)RZ(e)$.

\begin{lemma} \label{LZ_T,H_T,RZ_T lemma}
Let $S$ be a semigroup with idempotent $e$. If $T=LZ'(e)H'(e)RZ'(e)$ is a subsemigroup of $S$ for some $LZ'(e) \subset LZ(e)$ containing $e$, some subgroup $H'(e) \subset H(e)$, and some $RZ'(e) \subset RZ(e)$ containing $e$, then the maximum left zero subsemigroup of $T$ containing $e$ is $LZ_T(e)=LZ'(e)$, the maximum subgroup of $T$ containing $e$ is $H_T(e)=H'(e)$, and the maximum right zero subsemigroup of $T$ containing $e$ is $RZ_T(e)=RZ'(e)$.
\end{lemma}

\begin{proof}
Since $e$ is an idempotent contained in $LZ'(e)$, $H'(e)$, and $RZ'(e)$, $e=eee \in LZ'(e)H'(e)RZ'(e)=T$, and $T$ is a subsemigroup of $S$ containing $e$.

For any $e_L \in LZ'(e)$, $e_L=e_Lee \in T$, and $e_L \in LZ_T(e)$.
For any $e_L \in LZ'(e)$, $a \in H'(e)$, and $e_R \in RZ'(e)$, if $e_Lae_R \in LZ_T(e)$, $e(e_Lae_R)=e$. Since $ee_L=e$, $e(e_Lae_R)=(ee_L)ae_R=eae_R$, and $eae_R=e=eee$. By Lemma \ref{uniqueness lemma}, any product decomposition of $LZ(e)H(e)RZ(e)$ is unique, so $a=e=e_R$. Then $e_Lae_R=e_Lee=e_L \in LZ'(e)$. Therefore $LZ_T(e)=LZ'(e)$. Similarly, $RZ_T(e)=RZ'(e)$.

For any $a \in H'(e)$, there exists $a^{-1} \in H'(e)$. Since $e$ is the identity of $H'(e)$, $a=eae \in T$ and $a^{-1}=ea^{-1}e \in T$. Therefore $a \in H_T(e)$.
For any $e_L \in LZ'(e)$, $a \in H'(e)$, and $e_R \in RZ'(e)$, if $e_Lae_R \in H_T(e)$, then $e(e_Lae_R)=e_Lae_R$ and $(e_Lae_R)e=e_Lae_R$. Since $ee_L=e$ and $e_Re=e$, $eae_R=e_Lae_R$ and $e_Lae=e_Lae_R$. By Lemma \ref{uniqueness lemma}, $e=e_L$ and $e=e_R$.
Therefore $e_Lae_R=eae=a \in H'(e)$, and $H_T(e)=H'(e)$.
\end{proof}

\section{Right [Left] Subgroups}

A right [left] group is left [right] cancellative and right [left] simple semigroup. Let $S$ be a right group. For $a \in S$, there exists $b \in S$ such that $ab=a$ because $S$ is right simple. Then $abb=ab$, which implies $bb=b$ because $S$ is left cancellative. Therefore $S$ contains the idempotent $b$. Similarly, a left group contains an idempotent.

\begin{lemma} \label{right subgroup}
Let $S$ be a semigroup, and let $e \in S$ be an idempotent.
Let $T=H'(e)RZ'(e)$, for some subgroup $H'(e) \subset H(e)$ and some $RZ'(e) \subset RZ(e)$ with $e \in RZ'(e)$. Then $T$ is a right subgroup of $S$ containing $e$.
\end{lemma}

Similarly, $LZ'(e)H'(e)$, for some subgroup $H'(e) \subset H(e)$ and some $LZ'(e) \subset LZ(e)$ with $e \in LZ'(e)$, is a left subgroup of $S$ containing $e$.

Since $H'(e)$ is a subgroup of $H(e)$, $H'(e)$ contains $e$ as its identity. Since $e \in RZ'(e) \subset RZ(e)$, $RZ'(e)$ is a right zero subsemigroup of $S$ containing $e$.

\begin{proof}
Since $e$ is an idempotent, $e=ee \in H'(e)RZ'(e)=T$.
For any $x_1,x_2 \in H'(e)$ and $e_1,e_2 \in RZ'(e)$, 
$$(x_1e_1)(x_2e_2) = x_1e_1(ex_2)e_2 = x_1(e_1e)x_2e_2 = x_1ex_2e_2 = x_1x_2e_2 \in H'(e)RZ'(e).$$
Therefore $T$ is a subsemigroup of $S$ containing $e$.

For any $x_1,x_2,x_3 \in H'(e)$ and $e_1,e_2,e_3 \in RZ'(e)$,
if $(x_1e_1)(x_2e_2)=(x_1e_1)(x_3e_3)$, then $x_1x_2e_2=x_1x_3e_3$ as calculated above.
Since $x_1^{-1} \in H'(e)$, $x_1^{-1}x_1x_2e_2 = x_1^{-1}x_1x_3e_3$, which implies $x_2e_2 = x_3e_3$.
Therefore $T$ is left cancellative.

For any $x_1,x_2 \in H'(e)$ and $e_1,e_2 \in RZ'(e)$,
$$(x_1e_1)(x_1^{-1}x_2e_2)=x_1x_1^{-1}x_2e_2=ex_2e_2=x_2e_2.$$
Therefore $T$ is right simple, and thus it is a right group.
\end{proof}

\begin{theorem} \label{RG(e) maximum theorem}
Let $S$ be a semigroup, and let $e \in S$ be an idempotent. Then $RG(e)=H(e)RZ(e)$ is the maximum right subgroup of $S$ containing $e$.
\end{theorem}

Similarly, $LG(e)=LZ(e)H(e)$ is the maximum left subgroup of $S$ containing $e$.

\begin{proof}
By Lemma \ref{right subgroup}, $RG(e)$ is a right subgroup of $S$ containing $e$.

For any $x \in S$, if there exists $y \in S$ such that $x=ey$, then $ex=eey$, and since $e$ is an idempotent, $ex=eey=ey=x$, so that $e$ is a left identity of $x$.

Let $RG'(e)$ be a right subgroup of $S$ containing $e$. Since $e$ is an idempotent and right groups are right simple (and left cancellative), $e$ is a left identity of $RG(e)$ and $RG'(e)$. Let $a \in RG'(e)$. Since $a,e \in RG'(e)$, $ae \in RG'(e)$, and there exists $b \in RG'(e)$ such that $aeb=e$ because $RG'(e)$ is right simple. Therefore $(eae)(ebe)=e$. Similarly, there exists $c \in RG'(e)$ such that $bec=e$, so that $(ebe)(ece)=e$. Then $eae=(eae)(ebe)(ece)=ece$, and $eae \in H(e)$.
Since $eae \in H(e)$, there exists $(eae)^{-1}$, and $(eae)^{-1}ea \in RZ(e)$ because $e[(eae)^{-1}ea] = (eae)^{-1}ea$ and $[(eae)^{-1}ea]e=e$.
Thus $$a=ea=eea=(eae)[(eae)^{-1}ea] \in H(e)RZ(e)=RG(e),$$ and $RG(e)$ is the maximum right subgroup of $S$ containing $e$.
\end{proof}

If $S$ is a right [left] group, then it is its own maximal right [left] subsemigroup for every idempotent it contains.

\begin{corollary}
If $S$ is a right [left] group, then $S=H(e)RZ(e)$ [$S=LZ(e)H(e)$] for any idempotent $e \in S$.
\end{corollary}

We can then define another equivalence for right [left] groups.

\begin{corollary}
Let $S$ be a semigroup. There exists an idempotent $e \in S$ such that $S=H(e)RZ(e)$ [$S=LZ(e)H(e)$] if and only if $S$ is a right [left] group.
\end{corollary}

Let $S$ be a semigroup with idempotent $e$. Then for any subsemigroup $S'$ of $S$ with $S'$ containing $e$, $RZ'(e)$, the maximum right zero subsemigroup of $S'$ containing $e$, is a subset of $RZ(e)$ and contains $e$, and $H'(e)$, the maximum subgroup of $S'$ containing $e$, is a subgroup of $H(e)$. Then $H'(e)RZ'(e)=\{e\}H'(e)RZ'(e) \subset LZ(e)H(e)RZ(e)$, and uniqueness is implied by Lemma \ref{uniqueness lemma}.

\begin{corollary}
Let $S$ be a semigroup, let $e \in S$ be an idempotent.
Then $RG'(e)$ is a right subgroup of $S$ containing $e$ if and only if $RG'(e)=H'(e)RZ'(e)$ for some subgroup $H'(e) \subset H(e)$ and some $RZ'(e) \subset RZ(e)$ with $e \in RZ'(e)$. In this case, every product of $RG'(e)$ is unique, and $H'(e)$ and $RZ'(e)$ are unique.
\end{corollary}

Similarly, $LG'(e)$ is a left subgroup of $S$ containing $e$ if and only if $LG'(e)=LZ'(e)H'(e)$ for some unique subgroup $H'(e) \subset H(e)$ and some unique $LZ'(e) \subset LZ(e)$ with $e \in LZ'(e)$.

Let $E$ be the set of idempotents in semigroup $S$. Every right [left] subgroup contains an idempotent, and every idempotent in $S$ is contained in a maximal right [left] subgroup. Then $\{RG(e)|e \in E\}$ is the set of all maximal right subgroups in $S$, $\{LG(e)|e \in E\}$ is the set of all maximal left subgroups in $S$. For these sets and $e,f \in E$, either $RG(e)=RG(f)$ or $RG(e) \cap RG(f)= \emptyset$, and either $LG(e)=LG(f)$ or $LG(e) \cap LG(f)= \emptyset$.

\section{Rectangular Bands}

Every element of a rectangular band is an idempotent.

\begin{lemma} \label{rectangular band lemma}
Let $S$ be a semigroup, and let $e \in S$ be an idempotent.
Let $T=LZ'(e)RZ'(e)$ for some sets $LZ'(e)$ and $RZ'(e)$ such that
$e \in LZ'(e) \subset LZ(e)$, $e \in RZ'(e) \subset RZ(e)$, and $RZ'(e)LZ'(e) = \{e\}$.
Then $T$ is a rectangular band subsemigroup of $S$ containing $e$.
\end{lemma}

\begin{proof}
Since $e$ is an idempotent, $e=ee \in LZ'(e)RZ'(e)=T$.
For any $a,c \in LZ'(e)$ and $b,d \in RZ'(e)$, $ab,cd \in T$, and $$(ab)(cd)=a(bc)d=aed=(ae)d=ad \in T.$$ Therefore $T$ is a subsemigroup of $S$ containing $e$.

For any $a,c \in LZ'(e)$ and $b,d \in RZ'(e)$, $ab,cd \in T$, and $$(ab)(cd)(ab)=a(bc)(da)b=a(ee)b=aeb=a(eb)=ab.$$ Therefore $T$ is a rectangular band subsemigroup of $S$ containing $e$.
\end{proof}

\begin{theorem} \label{RB(e) theorem}
Let $S$ be a semigroup, and let $e \in S$ be an idempotent such that $RZ(e)LZ(e)=\{e\}$. Then $RB(e)=LZ(e)RZ(e)$ is the maximum rectangular band subsemigroup of $S$ containing $e$.
\end{theorem}

\begin{proof}
By Lemma \ref{rectangular band lemma}, $RB(e)=LZ(e)RZ(e)$ is a rectangular band subsemigroup of $S$ containing $e$. Let $U$ be a rectangular band subsemigroup of $S$ containing $e$. For any $a \in U$, $a=aea=aeea=(ae)(ea)$. Then $(ae)e=a(ee)=ae$ because $e$ is an idempotent, and $e(ae)=e$ because $a,e \in U$, which is a rectangular band. Therefore $ae \in LZ(e)$, and similarly $ea \in RZ(e)$. Thus $a \in LZ(e)RZ(e)=RB(e)$, and $RB(e)$ is the maximum rectangular band subsemigroup containing $e$.
\end{proof}

\begin{corollary} \label{S=RB(e) corollary}
Let $S$ be a rectangular band. For any $e \in S$, $LZ(e)=Se$, $RZ(e)=eS$, $RZ(e)LZ(e)=\{e\}$, and $S=RB(e)=LZ(e)RZ(e)$.
\end{corollary}

\begin{proof}
If $S$ is a rectangular band and $e \in S$, then $e$ is an idempotent. By Theorem \ref{left zero theorem} and its dual, $LZ(e)$ is the maximum left zero subsemigroup of $S$ containing $e$ and $RZ(e)$ is the maximum right zero subsemigroup of $S$ containing $e$.

For any $a \in S$, $(ae)e=a(ee)=ae$ and $e(ae)=e$. Therefore $ae \in LZ(e)$ and $Se \subset LZ(e)$. For any $b \in LZ(e)$, $b=be \in Se$, so $LZ(e) \subset Se$. Thus $LZ(e)=Se$. Similarly, $RZ(e)=eS$.

Since $S$ is a rectangular band, $RZ(e)LZ(e)=(eS)(Se)=e(SS)e=\{e\}$, and by Theorem \ref{RB(e) theorem}, $S=R(B)=LZ(e)RZ(e)=(Se)(eS)$.
\end{proof}

We can then define another equivalence for rectangular bands.

\begin{corollary} \label{rectangular band equiv corollary}
Let $S$ be a semigroup. There exists an idempotent $e \in S$ such that $S=LZ(e)RZ(e)$ if and only if $S$ is a rectangular band.
\end{corollary}

\begin{proof}
If $S=LZ(e)RZ(e)$ for idempotent $e \in S$, then for any $e_L \in LZ(e)$ and $e_R \in RZ(e)$, $e_Le=e_L \in S$ and $ee_R=e_R \in S$. Since $S$ is a semigroup, $e_Re_L \in S$ and $e_Re_L=f_Lf_R$ for some $f_L \in LZ(e)$ and $f_R \in RZ(e)$. Then $ee_Re_Le=ef_Lf_Re$, which simplifies to $e_Re_L=ee=e$. Thus $RZ(e)LZ(e)=\{e\}$, and $S$ is a rectangular band by Lemma \ref{rectangular band lemma}.

If $S$ is a rectangular band, then $S$ contains an idempotent $e$ and $S=LZ(e)RZ(e)$ by Corollary \ref{S=RB(e) corollary}.
\end{proof}

Let $S$ be a semigroup with idempotent $e$. Then for any rectangular band subsemigroup $S'$ of $S$ with $S'$ containing $e$, $RZ'(e)$, the maximum right zero subsemigroup of $S'$ containing $e$, is a subset of $RZ(e)$ and contains $e$, and $LZ'(e)$, the maximum left zero subsemigroup of $S'$ containing $e$, is a subset of $LZ(e)$ and contains $e$, such that $RZ'(e)LZ'(e)=\{e\}$. Then $LZ'(e)RZ'(e)=LZ'(e)\{e\}RZ'(e) \subset LZ(e)H(e)RZ(e)$, and uniqueness is implied by Lemma \ref{uniqueness lemma}.

\begin{corollary}
Let $S$ be a semigroup with idempotent $e$. $T$ is a rectangular band subsemigroup of $S$ containing $e$ if and only if $T=LZ'(e)RZ'(e)$ for some $LZ'(e) \subset LZ(e)$ with $e \in LZ'(e)$ and some $RZ'(e) \subset RZ(e)$ with $e \in RZ'(e)$. In this case, every product in $T$ is unique, and $LZ'(e)$ and $RZ'(e)$ are unique.
\end{corollary}

\section{Completely Simple Subsemigroups}

A completely simple semigroup is a simple semigroup containing a primitive idempotent.
There are several alternate definitions for a completely simple semigroup, including a completely simple semigroup is simple and contains a minimal left ideal and a minimal right ideal. There are also several well known results for completely simple semigroups, including every idempotent is a primitive idempotent and every principal left [right] ideal is a minimal left [right] ideal. See \cite{clifford},\cite{howie}, and \cite{grillet} for references.

\begin{lemma} \label{completely simple lemma}
Let $S$ be a semigroup, and let $e \in S$ be an idempotent.
Let $T=LZ'(e)H'(e)RZ'(e)$ for some $LZ'(e) \subset LZ(e)$ with $e \in LZ'(e)$, some subgroup $H'(e) \subset H(e)$, and some $RZ'(e) \subset RZ(e)$ with $e \in RZ'(e)$, such that $RZ'(e)LZ'(e) \subset H'(e)$. Then $T$ is a completely simple subsemigroup of $S$ containing $e$.
\end{lemma}

\begin{proof}
As an idempotent, $e=eee \in LZ'(e)H'(e)RZ'(e)=T$.

For any $e_L,f_L \in LZ'(e)$, $a,b \in H'(e)$, and $e_R,f_R \in RZ'(e)$,
$(e_Lae_R)(f_Lbf_R)=e_L(ae_Rf_Lb)f_R \in LZ'(e)H'(e)RZ'(e)=T$ because $e_Rf_L \in H'(e)$. Therefore $T$ is a subsemigroup of $S$. Also $f_Lbe,ea^{-1}f_R \in T$, and $(f_Lbe)(e_Lae_R)(ea^{-1}f_R)=f_L(beaea^{-1})f_R=f_Lbf_R$. Thus $T$ is simple.

For any idempotent $f \in T$ such that $f < e$, $ff=f$ and $fe=ef=f$. Also $f=e_Lae_R$ for some $a \in H'(e)$, and $ff=(e_Lae_R)(e_Lae_R)=e_Lae_R=f$. Therefore $e(e_Lae_R)(e_Lae_R)e=e(e_Lae_R)e$, which implies $ae_Re_La=a$, and $e_Re_L=a^{-1}$. Since $fe=f$ and $ef=f$, $(e_Lae_R)e=e_Lae_R$ and $e(e_Lae_R)=e_Lae_R$, or $e_La=e_Lae_R$ and $ae_R=e_Lae_R$. Therefore $a^{-1}e(e_La)=a^{-1}e(e_Lae_R)$ and $(ae_R)ea^{-1}=(e_Lae_R)ea^{-1}$, which implies $e=e_R$ and $e=e_L$. Thus $a^{-1}=ee=e$, and $f=eee=e$, so that $e$ is a primitive idempotent, and $T$ is completely simple.
\end{proof}

\begin{lemma} \label{eSe=H(e) lemma}
Let $S$ be a completely simple semigroup. Then for any idempotent $e \in S$, $eSe=H(e)$.
\end{lemma}

\begin{theorem} \label{CS(e) maximum theorem}
Let $S$ be a semigroup, and let $e \in S$ be an idempotent. If $RZ(e)LZ(e) \subset H(e)$, then $CS(e)=LZ(e)H(e)RZ(e)$ is the maximum completely simple subsemigroup of $S$ containing $e$.
\end{theorem}

\begin{proof}
By Lemma \ref{completely simple lemma}, $CS(e)$ is a completely simple subsemigroup of $S$ containing $e$.

Let $U$ be a completely simple subsemigroup of $S$ containing $e$, and let $a \in U$. Then by Lemma \ref{eSe=H(e) lemma}, $eae \in H(e)$, and $(eae)^{-1} \in H(e)$. Since $e(ae(eae)^{-1})=e$ and $(ae(eae)^{-1})e=ae(eae)^{-1}$, $ae(eae)^{-1} \in LZ(e)$, and similarly $(eae)^{-1}ea \in RZ(e)$.

Since $U$ is simple, there exist $x,y \in U$ such that $a=xey$. Then $eae=(exe)(eye)$, and again by Lemma \ref{eSe=H(e)^0 lemma}, $exe,eye \in H(e)$, so that $(eae)^{-1}=(eye)^{-1}(exe)^{-1}$. 

Consider $[ae(eae)^{-1}](eae)[(eae)^{-1}ea] \in LZ(e)H(e)RZ(e)$, which simplifies to $ae(eae)^{-1}ea$. Then 
$$ae(eae)^{-1}ea=xeye(eye)^{-1}(exe)^{-1}exey=xeey=xey=a.$$
Therefore $a \in CS(e)$.
\end{proof}

\begin{corollary} \label{completely simple corollary}
Let $S$ be a completely simple semigroup. Then for any idempotent $e \in S$, $RZ(e)LZ(e) \subset H(e)$ and $S=LZ(e)H(e)RZ(e)$.
\end{corollary}

\begin{proof}
For any $e_L \in LZ(e)$ and $e_R \in RZ(e)$, $e(e_Re_L)e=(ee_R)(e_Le)=e_Re_L \in eSe$, and $eSe=H(e)$ by Lemma \ref{eSe=H(e)^0 lemma}, so that $RZ(e)LZ(e) \subset H(e)$. By Theorem \ref{CS(e) maximum theorem}, $S=CS(e)=LZ(e)H(e)RZ(e)$.
\end{proof}

We can then define another equivalence for completely simple semigroups.

\begin{corollary} \label{completely simple equiv corollary}
Let $S$ be a semigroup. There exists an idempotent $e \in S$ such that $S=LZ(e)H(e)RZ(e)$ if and only if $S$ is completely simple.
\end{corollary}

\begin{proof}
Let $S$ be a semigroup with idempotent $e \in S$ such that $S=LZ(e)H(e)RZ(e)$.
For any $e_L \in LZ(e)$ and $e_R \in RZ(e)$, then $e_L=e_Lee \in S$ and $e_R=eee_R \in S$. Since $S$ is a semigroup, $e_Re_L \in S$, and $e_Re_L=f_Laf_R$ for some $f_L \in LZ(e)$, $a \in H(e)$, and $f_R \in RZ(e)$. Therefore 
$$e_Re_L=(ee_R)(e_Le)=e(e_Re_L)e=e(f_Laf_R)e=(ef_L)a(f_Re)=eae=a \in H(e).$$ 
By Lemma \ref{completely simple lemma}, $S$ is completely simple.

If $S$ is completely simple, then $S$ contains an idempotent $e$. By Corollary \ref{completely simple corollary}, $S=LZ(e)H(e)RZ(e)$.
\end{proof}

It is then apparent that a right zero semigroup $S$ is a completely simple semigroup such that for any idempotent $e \in S$, $LZ(e)=H(e)=\{e\}$. A right group $S$ is a completely simple semigroup such that for any idempotent $e \in S$, $LZ(e)=\{e\}$. A rectangular band $S$ is a completely simple semigroup such that for any idempotent $e \in S$, $H(e)=\{e\}$.

Lemma \ref{uniqueness lemma} applies to completely simple subsemigroups.

\begin{corollary}
Let $S$ be a semigroup with idempotent $e$. $T$ is a completely simple subsemigroup of $S$ containing $e$ if and only if $T=LZ'(e)H'(e)RZ'(e)$ for some $LZ'(e) \subset LZ(e)$ with $e \in LZ'(e)$, some subgroup $H'(e) \subset H(e)$, and some $RZ'(e) \subset RZ(e)$ with $e \in RZ'(e)$, such that $RZ'(e)LZ'(e) \subset H'(e)$. In this case, every product in $T$ is unique, and $LZ'(e)$, $H'(e)$, and $RZ'(e)$ are unique.
\end{corollary}

\section{Rectangular Subgroups}

A rectangular group is a semigroup that isomorphic to the direct product of a rectangular band and a group. Every rectangular band is a rectangular group using the trivial group. A rectangular band contains an idempotent $e$ and a group contains its identity $1$ as an idempotent. Therefore the direct product of a rectangular band and a group contains the idempotent $(e,1)$, and any rectangular group contains an idempotent.

\begin{lemma} \label{rectangular subgroup lemma}
Let $S$ be a semigroup, and let $e \in S$ be an idempotent.
Let $T=LZ'(e)H'(e)RZ'(e)$ for some $LZ'(e) \subset LZ(e)$ with $e \in LZ'(e)$, some subgroup $H'(e) \subset H(e)$, and some $RZ'(e) \subset RZ(e)$ with $e \in RZ'(e)$, such that $RZ'(e)LZ'(e)=\{e\}$. Then $T$ is a rectangular subgroup of $S$ containing $e$.
\end{lemma}

Since $e \in LZ'(e) \subset LZ(e)$, $LZ'(e)$ is a left zero subsemigroup of $S$ containing $e$. Similarly, $RZ'(e)$ is a right zero subsemigroup of $S$ containing $e$. Since $H'(e)$ is a subgroup of $H(e)$, $H'(e)$ contains $e$ as its identity. 

\begin{proof}
Let $T=LZ'(e)H'(e)RZ'(e)$. By Lemma \ref{completely simple lemma}, $T$ is a completely simple subsemigroup of $S$ containing $e$, and $U=LZ'(e)RZ'(e)$ is a rectangular band subsemigroup of $S$ and $T$ containing $e$ by Lemma \ref{rectangular band lemma}.

Define $F:T \to U \times H'(e)$ such that for any $e_Lae_R \in T$, $F(e_Lae_R)=(e_Le_R,a)$. By Lemma \ref{uniqueness lemma}, $F$ is well-defined and injective. It is also surjective. 

For any $e_Lae_R,f_Lbf_R \in T$, since $e_Rf_L \in RZ'(e)LZ'(e)=\{e\}$, 
$$(e_Lae_R)(f_Lbf_R)=e_La(e_Rf_L)bf_R=e_Laebf_R=e_Labf_R$$ and 
$$(e_Le_R,a)(f_Lf_R,b)=(e_Le_Rf_Lf_R,ab)=(e_Lef_R,ab)=(e_Lf_R,ab).$$ Therefore $F$ is an isomorphism, and $T$ is a rectangular group.
\end{proof}

\begin{lemma} \label{rectangular group lemma}
Let $S$ be a rectangular group. Then $S$ is completely simple.
\end{lemma}

\begin{proof}
As a rectangular group, $S$ is isomorphic to a direct product of a rectangular band and a group, $T$.
Let $T=U \times G$ for rectangular band $U$ and group $G$ with identity $1$. Since $U$ is a rectangular band, it contains an idempotent $e$, and $U=LZ_U(e)RZ_U(e)$ by Corollary \ref{rectangular band equiv corollary}. Then $(e,1)$ is an idempotent in $T$.

For any $(e_Le_R,a) \in T$, since $e_Re=e$ and every product decomposition in a rectangular band is unique, $(e_Le_R,a)(e,1)=(e_Le_R,a)$ if and only if $e=e_R$ and $(e_Le_R,a)=(e_L,a)$.
Also, $(e,1)(e_L,a)=(e,1)$ if and only if $a=1$. Therefore, $LZ_T(e,1)=LZ_U(e) \times \{1\}$. Similarly, $RZ_T(e,1)=RZ_U(e) \times \{1\}$.

For any $(e_Le_R,a) \in T$, $(e_Le_R,a)(e,1)=(e_Le_R,a)$ and $(e,1)(e_Le_R,a)=(e_Le_R,a)$ if and only if $e_L=e_R=e$, so that $(e_Le_R,a)=(e,a)$. Then $(e,a)(e,a^{-1})=(e,a^{-1})(e,a)=(e,1)$, which implies $(e,a) \in H_T(e,1)$. Therefore $H_T(e,1)=\{e\} \times G$.

For any $(e_Le_R,a) \in T$, $(e_Le_R,a)=(e_L,1)(e,a)(e_R,1) \in LZ_T(e,1)H_T(e,1)RZ_T(e,1)$. Thus $T=LZ_T(e,1)H_T(e,1)RZ_T(e,1)$, and by Corollary \ref{completely simple equiv corollary}, $T$ is completely simple, which implies $S$ is completely simple.
\end{proof}

\begin{theorem} \label{RECTG(e) maximum theorem}
Let $S$ be a semigroup, and let $e \in S$ be an idempotent such that $RZ(e)LZ(e)=\{e\}$.
Then $RECTG(e)=LZ(e)H(e)RZ(e)$ is the maximum rectangular subgroup of $S$ containing $e$.
\end{theorem}

\begin{proof}
By Lemma \ref{rectangular subgroup lemma}, $RECTG(e)$ is a rectangular subgroup of $S$ containing $e$, and by Theorem \ref{CS(e) maximum theorem}, $RECTG(e)=CS(e)$ is the maximum completely simple subsemigroup of $S$ containing $e$.

For any rectangular subgroup of $S$ containing $e$, it is completely simple by Lemma \ref{rectangular group lemma} and contained in $RECTG(e)$.
Therefore $RECTG(e)$ is the maximum rectangular subgroup of $S$ containing $e$.
\end{proof}

\begin{corollary}
Let $S$ be a rectangular group. Then for any idempotent $e \in S$, $RZ(e)LZ(e)=\{e\}$ and $S=LZ(e)H(e)RZ(e)$.
\end{corollary}

\begin{proof}
If $S$ is a rectangular group, then $S$ is isomorphic to $T=U \times G$ for rectangular band $U$ and group $G$ with identity $1$. Any idempotent in $T$ is of the form $(f,1)$ for some idempotent $f \in U$. From the proof of Lemma \ref{rectangular group lemma}, $LZ_T(f,1)=LZ_U(f) \times \{1\}$ and $RZ_T(f,1)=RZ_U(f) \times \{1\}$. Therefore for any $f_L \in LZ_U(f)$ and $f_R \in RZ_U(f)$, $f_Rf_L=f$ because $U$ is a rectangular band, so that $RZ_T(f,1)LZ_T(f,1)=\{(f,1)\}$. 

Since $S$ is isomorphic to $T$, for any idempotent $e \in S$, $RZ(e)LZ(e)=\{e\}$, and by Theorem \ref{RECTG(e) maximum theorem}, $S=LZ(e)H(e)RZ(e)$.
\end{proof}

We can then define another equivalence for rectangular groups.

\begin{corollary}
Let $S$ be a semigroup. There exists an idempotent $e \in S$ such that $S=LZ(e)H(e)RZ(e)$ and $RZ(e)LZ(e)=\{e\}$ if and only if $S$ is a rectangular group.
\end{corollary}

Lemma \ref{uniqueness lemma} applies to rectangular subgroups.

\begin{corollary}
Let $S$ be a semigroup with idempotent $e$. $T$ is a rectangular subgroup of $S$ containing $e$ if and only if $T=LZ'(e)H'(e)RZ'(e)$ for some $LZ'(e) \subset LZ(e)$ with $e \in LZ'(e)$, some subgroup $H'(e) \subset H(e)$, and some $RZ'(e) \subset RZ(e)$ with $e \in RZ'(e)$, such that $RZ'(e)LZ'(e)=\{e\}$. In this case, every product in $T$ is unique, and $LZ'(e)$, $H'(e)$, and $RZ'(e)$ are unique.
\end{corollary}

\section{Completely $0$-Simple Subsemigroups}

A completely $0$-simple semigroup is $0$-simple and contains a primitive idempotent. There are several alternate definitions, including a completely $0$-simple semigroup is $0$-simple and contains a $0$-minimal left ideal and a $0$-minimal right ideal. There are also several well known results for completely $0$-simple semigroups, including every nonzero idempotent is a primitive idempotent and every principal left [right] ideal generated by a nonzero element is a $0$-minimal left [right] ideal. See \cite{clifford},\cite{howie}, and \cite{grillet} for references.

Define $L(e)=\{x \in L_e|\text{$ex=e$ or $ex=0$}\}$. Since $x \in L_e$, $xe=x$, so $L(e)=LZ(e) \cup \{x \in L_e|ex=0\}$. Let $lzediv(e)=\{x \in L_e|ex=0\}$. Then $x \in lzediv(e)$ if and only if $xe=x$, $ex=0$, and there exists $w \in S$ such that $wx=e$. Similarly, define $R(e)=\{y \in R_e|\text{$ye=e$ or $ye=0$}\}=RZ(e) \cup rzediv(e)$, and $y \in rzediv(e)$ if and only if $ey=y$, $ye=0$, and there exists $z \in S$ such that $yz=e$.

\begin{lemma} \label{completely $0$-simple lemma}
Let $S$ be a semigroup with zero, and let $e \in S$ be a nonzero idempotent.
Let $T=L'(e)H'(e)^0R'(e)$ for some $L'(e) \subset L(e)$ with $e \in L'(e)$, some $H'(e)$ subgroup of $H(e)$, and some $R'(e) \subset R(e)$ with $e \in R'(e)$, such that $R'(e)L'(e) \subset H'(e)^0$, $R'(e)x \neq \{0\}$ for all $x \in L'(e)$, and $yL'(e) \neq \{0\}$ for all $y \in R'(e)$.
Then $T$ is a completely $0$-simple subsemigroup of $S$ containing $e$.
\end{lemma}

Since $e \in L'(e),R'(e)$, for any $x \in L'(e) \cap LZ(e)$, $ex=e \neq 0$, and for any $y \in R'(e) \cap RZ(e)$, $ye=e \neq 0$, so it is necessarily true that $R'(e)x \neq \{0\}$ for all $x \in L'(e) \cap LZ(e)$ and $yL'(e) \neq \{0\}$ for all $y \in L'(e) \cap RZ(e)$.

\begin{proof}
$e=eee \in L'(e)H'(e)^0R'(e)=T$.

For any $x,s \in L'(e)$, $a,b \in H'(e)^0$, and $y,t \in R'(e)$,
$(xay)(sbt)=xa(ys)bt \in L'(e)H'(e)^0R'(e)=T$ because $ys \in H'(e)^0$.
Therefore $T$ is a subsemigroup of $S$.

Assume $xay \neq 0$, then $a \neq 0$ and $a^{-1} \in H'(e)$. 
There exists $w \in R'(e)$ such that $wx \neq 0$ and $wx \in H'(e)$ because $R'(e)L'(e) \subset H'(e)^0$. Let $X=sb(wx)^{-1}w$. There exists $z \in L'(e)$ such that $yz \neq 0$ and $yz \in H'(e)$ because $R'(e)L'(e) \subset H'(e)^0$. Let $Y=z(yz)^{-1}a^{-1}t$. Then 
$$X(xay)Y = (sb(wx)^{-1}w)xay(z(yz)^{-1}a^{-1}t)=sbeaea^{-1}t = sbt.$$
Thus $T$ is $0$-simple.

For any nonzero idempotent $f \in T$ such that $f < e$, $ff=f$ and $fe=ef=f$. Let $f=xay$ with nonzero $a \in H'(e)$. Then $xayxay=xay$ implies $yx \in H'(e)$ because $f$ is nonzero and $R'(e)L'(e) \subset H'(e)^0$. Since $xayxay=xay$, then $yxayxayx=yxayx$, which implies $yx=a^{-1}$. Since $ef=f$, $exay=xay$, and therefore $x \in LZ(e)$ and $ay=xay$. Similarly, since $fe=f$, $xaye=xay$, and therefore $y \in RZ(e)$ and $xa=xay$. By the uniqueness of product decomposition in $LZ(e)H(e)RZ(e)$, $x=y=e$, so that $a^{-1}=e$ and $a=e$. Then $f=eee=e$, and $e$ is a primitive idempotent, so that $T$ is completely $0$-simple.
\end{proof}

\begin{lemma} \label{eSe=H(e)^0 lemma}
Let $S$ be a completely $0$-simple semigroup. Then for any nonzero idempotent $e \in S$, $eSe=H(e)^0$.
\end{lemma}

\begin{theorem} \label{CZS(e) maximum theorem}
Let $S$ be a semigroup with zero, and let $e \in S$ be a nonzero idempotent.
If $R(e)L(e) \subset H(e)^0$, $R(e)x \neq \{0\}$ for any $x \in L(e)$, and $yL(e) \neq \{0\}$ for any $y \in R(e)$, then $CZS(e)=L(e)H(e)^0R(e)$ is the maximum completely $0$-simple subsemigroup of $S$ containing $e$.
\end{theorem}

\begin{proof}
By Lemma \ref{completely $0$-simple lemma}, $CZS(e)$ is a completely $0$-simple subsemigroup of $S$ containing $e$.

Let $U$ be a completely $0$-simple subsemigroup of $S$ containing $e$, and let $a \in U$ be nonzero. Since $U$ is $0$-simple, there exist $x,y \in U$ such that $a=xey$. 
Since $a \neq 0$, $a=xeey$ implies $xe,ey \neq 0$, and since $e$ is an idempotent $(xe)e=xe$ and $e(ey)=ey$. Since $xe \neq 0$, $Ue=Uxe$ because $U$ is completely $0$-simple, and there exists $w \in U$ with $wxe=e$. Therefore $xe \in L_e$. 
If $exe \neq 0$, then $exe \in H(e)$ by Lemma \ref{eSe=H(e)^0 lemma}, and if $exe=0$, then $xe \in lzediv(e)$.  
Similarly, if $eye \neq 0$, then $eye \in H(e)$, and if $eye=0$, then $ye \in rzediv(e)$.

\begin{itemize}
\item[] If $exe,eye \neq 0$, then
\begin{align*}
a &= (xe)ee(ey) \\
&= [xe(exe)^{-1}][(exe)(eye)][(eye)^{-1}ey] \\
&\in LZ(e)H(e)RZ(e) \subset CZS(e).
\end{align*}
\item[] If $exe=eye=0$, then
\begin{align*}
a &= (xe)e(ey) \\
&\in lzediv(e)H(e)rzediv(e) \subset CZS(e).
\end{align*}
\item[] If $exe=0$ and $eye \neq 0$, then 
\begin{align*}
a &= (xe)e(ey) \\
&= [xe][(eye)][(eye)^{-1}ey] \\
&\in lzediv(e)H(e)RZ(e) \subset CZS(e).
\end{align*}
\item[] If $exe \neq 0$ and $eye=0$, then
\begin{align*}
a &= (xe)e(ey) \\
&= [xe(exe)^{-1}][(exe)][ey] \\
&\in LZ(e)H(e)rzediv(e) \subset CZS(e).
\end{align*}
\end{itemize}

Therefore $CZS(e)$ is the maximum completely $0$-simple subsemigroup of $S$ containing $e$.
\end{proof}

\begin{corollary} \label{completely $0$-simple corollary}
Let $S$ be a completely $0$-simple semigroup. Then for any nonzero idempotent $e \in S$, $R(e)L(e) \subset H(e)^0$, $R(e)x \neq \{0\}$ for any $x \in L(e)$, $yL(e) \neq \{0\}$ for any $y \in R(e)$, and $S=L(e)H(e)^0R(e)$.
\end{corollary}

\begin{proof}
For any $x \in L(e)$ and $y \in R(e)$, $xe=x$ and $ey=y$, so $yx=eyxe \in eSe=H(e)^0$ by Lemma \ref{eSe=H(e)^0 lemma}.

Since $e \in LZ(e),RZ(e)$, if $x \in LZ(e)$, then $ex=e \neq 0$ and $R(e)x \neq \{0\}$, and if $y \in RZ(e)$, then $ye=e \neq 0$ and $yL(e) \neq \{0\}$.

Now assume $x \in lzediv(e)$ and $y \in rzediv(e)$. Since $x \in L_e$, there exists $w \in S$ such that $wx=e$. Then $ewx=e$, $eew=ew$, and $ew \in R_e$. Now either $ewe \neq 0$ or $ewe=0$. If $ewe \neq 0$, then $ewe \in H(e)$ by Lemma \ref{eSe=H(e)^0 lemma}, and $(ewe)^{-1}ew \in RZ(e) \subset R(e)$, and $(ewe)^{-1}ewx=(ewe)^{-1} \neq 0$. If $ewe=0$, then $ew \in rzediv(e) \subset R(e)$, and $ewx=e \neq 0$. Therefore $R(e)x \neq \{0\}$. Similarly, $yL(e) \neq \{0\}$.

By Theorem \ref{CZS(e) maximum theorem}, $S=CZS(e)=L(e)H(e)^0R(e)$. 
\end{proof}

We can then define another equivalence for completely $0$-simple semigroups.

\begin{corollary} \label{completely $0$-simple equiv corollary}
Let $S$ be a semigroup with $0$. There exists a nonzero idempotent $e \in S$ such that $S=L(e)H(e)^0R(e)$ if and only if $S$ is completely $0$-simple.
\end{corollary}

\begin{proof}
Let $S$ be a semigroup with $0$ and with a nonzero idempotent $e \in S$ such that $S=L(e)H(e)^0R(e)$.
For any $x \in L(e)$ and $y \in R(e)$, $x=xee \in S$ and $y=eey \in S$. Since $S$ is a semigroup, $yx \in S$, and $yx=ras$ for some $r \in L(e)$, $a \in H(e)^0$, and $s \in R(e)$ because $S$ is a semigroup. If $yx=0$, then $yx \in H(e)^0$. Consider $yx \neq 0$. Then $eyxe=erase$, which implies $yx=erase$ with $er \neq 0$ and $se \neq 0$. Therefore $r \in LZ(e)$ and $s \in RZ(e)$, so that $er=e$ and $se=e$, and $yx=eae=a \in H(e)$. 

Also for $x \in L(e)$, if $x \in LZ(e)$, then $ex=e \neq 0$, and $R(e)x \neq \{0\}$. If $x \in lzediv(e)$, then for some $t \in L(e)$, $b \in H(e)^0$, and $u \in R(e)$, $tbux=e$ because $x \in L_e$. Then $ux \neq 0$, so that $R(e)x \neq \{0\}$. Similarly, for $y \in R(e)$, $yL(e) \neq \{0\}$, and by Lemma \ref{completely $0$-simple lemma}, $S$ is completely $0$-simple.

If $S$ is completely $0$-simple, then $S$ contains a nonzero idempotent $e$. By Corollary \ref{completely $0$-simple corollary}, $S=L(e)H(e)^0R(e)$.
\end{proof}

\begin{corollary}
Let $S$ be a semigroup with $0$ and with nonzero idempotent $e$. $T$ is a completely $0$-simple subsemigroup of $S$ containing $e$ if and only if $T=L'(e)H'(e)^0R'(e)$
for some $L'(e) \subset L(e)$ with $e \in L'(e)$, some $H'(e)$ subgroup of $H(e)$, some $R'(e) \subset R(e)$ with $e \in R'(e)$, such that $R'(e)L'(e) \subset H'(e)^0$, $R'(e)x \neq \{0\}$ for all $x \in L'(e)$, and $yL'(e) \neq \{0\}$ for all $y \in R'(e)$.
\end{corollary}

Here $LZ'(e)=L'(e) \cap LZ(e)$, $H'(e)$, and $RZ'(e)=R'(e) \cap RZ(e)$ are unique, however $lzediv'(e)=L'(e) \cap lzediv(e)$ and $rzediv'(e)=R'(e) \cap rzediv(e)$ are not necessarily unique. We can make $lzediv'(e)$ and $rzediv'(e)$ unique if we require $lzediv'(e)H'(e)=lzediv'(e)$ and $H'(e)rzediv'(e)=rzediv'(e)$, but products in $L'(e)H'(e)^0R'(e)$ will not necessarily be unique.

For semigroup $S$ with element $a$, let $H_a=\{b \in S|S^1a=S^1b, \, aS^1=bS^1\}$ be the $\mathscr{H}$-class of $a$.

\begin{lemma} \label{H_a lemma}
Let $S$ be a semigroup with idempotent $e$. For any $a \in L_e$, $H_a = aH(e)$, and for any $a \in R_e$, $H_a = H(e)a$.
\end{lemma}

The above lemma is a result of Green's Lemma.

\begin{lemma} \label{L_e, R_e lemma}
Let $S$ be a completely $0$-simple semigroup with idempotent $e$. Then $L_e=LZ(e)H(e) \cup lzediv(e)$ and $R_e=H(e)RZ(e) \cup rzediv(e)$.
\end{lemma}

\begin{proof}
If $x \in L_e$, then $xe=x$ because $e$ is an idempotent.
Since $S$ is completely $0$-simple, $exe \in H(e)^0$ by Lemma \ref{eSe=H(e)^0 lemma}. If $exe=0$, then $x \in lzediv(e)$. If $exe \neq 0$, then $exe \in H(e)$, and $xe(exe)^{-1} \in LZ(e)$, which implies
\begin{align*}
x &= xee \\
&= [xe(exe)^{-1}][(exe)] \\
&\in LZ(e)H(e).
\end{align*}
Thus $L_e \subset LZ(e)H(e) \cup lzediv(e)$.

If $f \in LZ(e)$, $ef=e$ and $fe=f$, and $f \in L_e$. Then $H_f=fH(e)$ by Lemma \ref{H_a lemma}, and $LZ(e)H(e) \subset L_e$. By definition, $lzediv(e) \subset L_e$. Therefore $LZ(e)H(e) \cup lzediv(e) \subset L_e$.
\end{proof}

The following is equivalent to the Rees-Sushkevich Theorem.

\begin{theorem}[Rees-Sushkevich Theorem]
Let $S$ be a completely $0$-simple semigroup. For any nonzero idempotent $e \in S$,
let $I$ be an index of $\mathscr{R}$-classes in $D_e$ and let $L'(e)=\{x_i|i \in I\}$ such that $x_i \in L_e \cap R_i$ for any $i \in I$. Let $\Lambda$ be an index of $\mathscr{L}$-classes in $D_e$ and let $R'(e)=\{y_{\lambda}|\lambda \in \Lambda\}$ such that $y_{\lambda} \in R_e \cap L_{\lambda}$ for any $\lambda \in \Lambda$.
Then $R'(e)L'(e) \subset H(e)^0$, $R'(e)x \neq \{0\}$ for any $x \in L'(e)$, $yL'(e) \neq \{0\}$ for any $y \in R'(e)$, and $S=L'(e)H(e)^0R'(e)$, and every nonzero element of $S$ has a unique product decomposition in $L'(e)$, $H(e)$, and $R'(e)$.
\end{theorem}

\section{Rectangular $0$-Bands}

\begin{definition}
A {\bf rectangular $0$-band} is a completely $0$-simple semigroup $S$ such that for any $a,b \in S$, $aba=a$ or $aba=0$.
\end{definition}

Since it is completely $0$-simple, a rectangular $0$-band contains $0$ and a nonzero idempotent. The smallest rectangular $0$-band is the adjunction of the trivial group with $0$. For any semigroup $S$ with nonzero idempotent $e$, the smallest rectangular $0$-band subsemigroup containing $e$ is $\{e\}^0$.

There are other statements equivalent to this definition of rectangular $0$-band.

\begin{theorem} \label{rectangular $0$-band theorem}
Let $S$ be a semigroup with zero. The following are equivalent.
\begin{itemize}
\item[(1)] $S$ is completely $0$-simple and for any $a,b \in S$, $aba=a$ or $aba=0$.
\item[(2)] $S$ is completely $0$-simple and for any $a \in S$, $a^2=a$ or $a^2=0$.
\item[(3)] There exists a nonzero idempotent $e \in S$ such that $S=L(e)\{e\}^0R(e)$.
\item[(4)] There exists a nonzero idempotent in $S$ and for any nonzero idempotent $e \in S$, $S=L(e)\{e\}^0R(e)$.
\end{itemize}
\end{theorem}

\begin{proof}
By Corollary \ref{completely $0$-simple equiv corollary}, if $S$ is completely $0$-simple, then there exists a nonzero idempotent $e \in S$ such that $S=L(e)H(e)^0R(e)$. Let $a \in H(e)$. If $(1)$ is given, then $aea=a$ or $aea=0$, which implies $a^2=a$ and $a=e$, so that $H(e)=\{e\}$. If $(2)$ is given, then $a^2=a$ or $a^2=0$, which implies $a^2=a$ and $a=e$, so that $H(e)=\{e\}$. Therefore $(1)$ and $(2)$ imply $(3)$.

If there exists a nonzero idempotent $e \in S$ such that $S=L(e)\{e\}^0R(e)$, then for any $a \in H(e)$, $a=xey=xy$ for some $x \in L(e)$ and $y \in R(e)$. Since $ea=a=ae$, $exy=xy=xye$, and $ex,ye \neq 0$, so that $x \in LZ(e)$ and $y \in RZ(e)$. Therefore $a=eae=exye=ee=e$. Thus $H(e)=\{e\}$, and $S$ is a compeletely $0$-simple semigroup by Corollary \ref{completely $0$-simple equiv corollary}, and $R(e)L(e) \subset \{e\}^0$ by Corollary \ref{completely $0$-simple corollary}.
For any $x,s \in L(e)$ and $y,t \in R(e)$, either $ys=e$ or $ys=0$. If $ys=e$, then $(xey)(set)=xet$, and if $ys=0$, then $(xey)(set)=0$. Therefore $(xey)(set)(xey)=xey$ or $(xey)(set)(xey)=0$, and $(xey)(xey)=xey$ or $(xey)(xey)=0$. Thus $(3)$ implies $(1)$ and $(2)$.

Also, given $(3)$ which also implies $(2)$, $S$ is completely $0$-simple, and by Corollary \ref{completely $0$-simple corollary}, for any nonzero idempotent $f \in S$, $S=L(f)H(f)^0R(f)$. Then because $a^2=a$ or $a^2=0$ for any $a \in H(f)$, $a=f$ and $H(f)=\{f\}$, which proves $(4)$.

It is trivially true that $(4)$ implies $(3)$.
\end{proof}

A rectangular $0$-band $S$ is a completely $0$-simple semigroup such that for any idempotent $e \in S$, $H(e)=\{e\}$.

\begin{corollary} \label{rectangular $0$-band subsemigroup corollary}
Let $S$ be a semigroup with $0$ and with nonzero idempotent $e$. $T$ is a rectangular $0$-band subsemigroup of $S$ containing $e$ if and only if $T=L'(e)\{e\}^0R'(e)$
for some $L'(e) \subset L(e)$ with $e \in L'(e)$ and some $R'(e) \subset R(e)$ with $e \in R'(e)$, such that $R'(e)L'(e) \subset \{e\}^0$, $R'(e)x \neq \{0\}$ for all $x \in L'(e)$, and $yL'(e) \neq \{0\}$ for all $y \in R'(e)$. Every nonzero element of $T$ has a unique product decomposition in $L'(e)$, $\{e\}$, and $R'(e)$, and $L'(e)=L_T(e)$ and $R'(e)=R_T(e)$ are unique to $T$.
\end{corollary}

\begin{proof}
We only prove the uniqueness of product decomposition because the rest of the corollary follows directly from previous results.

Let $x,s \in L'(e)$ and $y,t \in R'(e)$ such that $xey=set$. Then there exists $z \in R'(e)$ such that $yz \neq 0$, which implies $yz=e$ because $R'(e)L'(e) \subset \{e\}^0$. Then $xeyz=setz$, or $x=setz$.
Therefore $tz \neq 0$, and $tz=e$ because $R'(e)L'(e) \subset \{e\}^0$, so that $x=s$. Similarly, $y=t$.
\end{proof}

\section{Rectangular $0$-Groups}

\begin{definition}
A {\bf rectangular $0$-group} is a semigroup which is isomorphic to the quotient of a direct product of a rectangular $0$-band and a group by its kernel.
\end{definition}

In this case, the kernel is the direct product of $\{0\}$ and the group.
Every rectangular $0$-band is a rectangular $0$-group using the trivial group.
A rectangular $0$-band contains a nonzero idempotent $e$ and a group contains its identity $1$ as an idempotent. Therefore the quotient of the direct product of a rectangular $0$-band and a group contains the nonzero idempotent $(e,1)$, and any rectangular $0$-group contains a nonzero idempotent. 

There are other statements equivalent to the definition of rectangular $0$-group.

\begin{theorem}
Let $S$ be a semigroup with zero. The following are equivalent.
\begin{itemize}
\item[(1)] $S$ is rectangular $0$-group.
\item[(2)] There exists a nonzero idempotent in $S$ and for any nonzero idempotent $e \in S$, $S=L'(e)H(e)^0R'(e)$ for some $L'(e)=LZ(e) \cup lzediv'(e)$ with $lzediv'(e) \subset lzediv(e)$ and some $R'(e)=RZ(e) \cup rzediv'(e)$ with $rzediv'(e) \subset rzediv(e)$ such that $R'(e)L'(e) \subset \{e\}^0$.
\item[(3)] There exists a nonzero idempotent $e \in S$ such that $S=L'(e)H(e)^0R'(e)$ for some $L'(e)=LZ(e) \cup lzediv'(e)$ with $lzediv'(e) \subset lzediv(e)$ and some $R'(e)=RZ(e) \cup rzediv'(e)$ with $rzediv'(e) \subset rzediv(e)$ such that $R'(e)L'(e) \subset \{e\}^0$.
\end{itemize}
\end{theorem}

\begin{proof}
$(1) \implies (2)$: As a rectangular $0$-group, $S$ is isomorphic to $T$, the quotient of a direct product of a rectangular $0$-band and a group by its kernel.
Let $T=(U \times G)/(\{0\} \times G)$ for rectangular $0$-band $U$ and group $G$. For any nonzero idempotent $(x,a) \in T$, $(x,a)(x,a)=(x,a)$ if and only if $xx=x \neq 0$ and $aa=a$, so $(x,a)$ is a nonzero idempotent if and only if $x$ is a nonzero idempotent and $a=1$. Therefore $E_T=E_U \times \{1\}$.

By Theorem \ref{rectangular $0$-band theorem}, as a rectangular $0$-band, $U$ contains a nonzero idempotent and for any nonzero idempotent $e \in U$, $U=L_U(e)\{e\}^0R_U(e)$. 
Then $T$ contains a nonzero idempotent and every nonzero idempotent $T$ is of the form $(e,1)$.

For any nonzero $x \in U$, $(x,a)(e,1)=(xe,a)=(x,a)$ if and only if $xe=x$. For any nonzero $x \in U$, $(e,1)(x,a)=(ex,a)=(e,1)$ if and only if $ex=e$ and $a=1$, and for any nonzero $x \in U$, $(e,1)(x,a)=(ex,a)=(0,1)$ if and only if $ex=0$. Also, there exists $(y,b) \in U \times G$ such that $(y,b)(x,a)=(e,1)$ if and only if $yx=e$ and $b=a^{-1}$.

Therefore $LZ_T(e,1)=LZ_U(e) \times \{1\}$ and $lzediv_U(e) \times \{1\} \subset lzediv_T(e,1)=lzediv_U(e) \times G$, and we let $lzediv_T'(e,1)=lzediv_U(e) \times \{1\}$. Similarly, $RZ_T(e,1)=RZ_U(e) \times \{1\}$ and $rzediv_U(e) \times \{1\} \subset rzediv_T(e,1)=rzediv_U(e) \times G$, and we let $rzediv_T'(e,1)=rzediv_U(e) \times \{1\}$.

For any $(x,1) \in L_T'(e,1)$ and $(y,1) \in R_T'(e,1)$, $(y,1)(x,1)=(yx,1)$. Since $R_U(e)L_U(e) \subset \{e\}^0$, either $yx=e$ and $(yx,1)=(e,1)$ or $yx=0$ and $(yx,1)=(0,1)$. Therefore $R_T'(e,1)L_T'(e,1) \subset \{(e,1)\}^0$.

For any nonzero $(xy,a) \in T$, $(xy,a)(e,1)=(xy,a)$ and $(e,1)(xy,a)=(xy,a)$ if and only if $x=y=e$, so that $(xy,a)=(e,a)$. Then $(e,a)(e,a^{-1})=(e,a^{-1})(e,a)=(e,1)$, which implies $(e,a) \in H_T(e,1)$. Therefore $H_T(e,1)=\{e\} \times G$.

For any nonzero $(xy,a) \in T$, $(xy,a)=(x,1)(e,a)(y,1) \in L_T'(e,1)H_T(e,1)R_T'(e,1)$, so that $T=L_T'(e,1)H_T(e,1)^0R_T'(e,1)$, and similar results apply to $S$.

$(2) \implies (3)$: This is trivially true.

$(3) \implies (1)$: Let $U=L'(e)\{e\}^0R'(e)$. Then $U \subset S$, and $e=eee \in L'(e)\{e\}^0R'(e) = U$. For any $x,s \in L'(e)$, and $y,t \in R'(e)$,
$(xey)(set)=(xy)(st)=x(ys)t \in L'(e)\{e\}^0R'(e)=U$ because $ys \in \{e\}^0$.
Therefore $U$ is a subsemigroup of $S$ containing $e$.

For any nonzero $xey=xy \in L'(e)\{e\}^0R'(e)=U$, $(xy)e=xy$ if and only if $y=e$ and $xy=x$. Then $ex=e$ if and only if $x \in LZ(e)$, so that $LZ_U(e)=LZ(e)$. 
For $x \in lzediv'(e)$, $ex=0$. Then $sbt \in L'(e)H(e)^0R'(e)=S$ such that $sbtx=e$ if and only if $tx=e$ and $s=b=e$. For $x \in lzediv'(e)$, there exists $t \in R'(e)$ such that $tx=e$, which implies $lzediv_U(e)=lzediv'(e)$. Therefore $L_U(e)=LZ_U(e) \cup lzediv_U(e)=L'(e)$. Similarly $R_U(e)=R'(e)$.
Also for nonzero $xey=xy \in L'(e)\{e\}^0R'(e)=U$, $(xy)e=xy$ and $e(xy)=xy$ if and only if $x=y=e$ and $xy=ee=e$, so that $H_U(e)=\{e\}$. Thus $U=L_U(e)\{e\}^0R_U(e)$, and $U$ is a rectangular $0$-band by Theorem \ref{rectangular $0$-band theorem}.

Define a function $F:S \to (U \times H(e))/(\{0\} \times H(e))$ such that for any $xay \in L'(e)H(e)^0R'(e)$, $F(xay)=(xy,a)$ if $a \neq 0$ and $F(xay)=(0,e)$ if $a = 0$. Let $xay,sbt \in L'(e)H(e)^0R'(e)$. There exists $w \in S$ such that $wx=e$ because $x \in L_e$, and there exists $z \in S$ such that $yz=e$ because $y \in R_e$.

Let $xay,sbt \neq 0$. If $xay=sbt$, then $wxayz=wsbtz$, which simplifies to $eae=a=wsbtz$. Since $a \neq 0$, $ws$ and $tz$ are nonzero and must equal $e$ because $R'(e)L'(e) \subset \{e\}^0$. Therefore $a=ebe=b$. Thus $a^{-1}wxay=a^{-1}wsat$ and $xayza^{-1}=satza^{-1}$, so that $y=t$ and $x=s$. Therefore every nonzero product in $L'(e)H(e)^0R'(e)$ is unique, and the function $F$ is well-defined and injective.
For any nonzero $(xy,a) \in (U \times H(e))/(\{0\} \times H(e))$, $F(xay)=(xy,a)$, so that $F$ is surjective.

Since $R'(e)L'(e) \subset \{e\}^0$, $ys \in \{e\}^0$, and $xaysbt=xaebt=xabt$ when $ys=e$ and $xaysbt=0$ when $ys=0$. 
If $ys=e$ and $ab \neq 0$, then 
$F(xaysbt)=F(xabt)=(xt,ab)$, and $F(xay)F(sbt)=(xy,a)(st,b)=(xyst,ab)=(xt,ab)$.
If $ys=0$ or $ab=0$, then $F(xaysbt)=(0,e)=F(xay)F(sbt)$.
Therefore $F$ is an isomorphism, and $S$ is a rectangular $0$-group.
\end{proof}

A rectangular $0$-group $S$ is a completely $0$-simple semigroup such that for any idempotent $e \in S$, there exist $L'(e) \subset L(e)$ and $R'(e) \subset R(e)$ such that $L'(e)H(e)=L(e)$, $H(e)R'(e)=R(e)$, and $R'(e)L'(e) \subset \{e\}^0$.

\begin{corollary} \label{rectangular $0$-subgroup corollary}
Let $S$ be a semigroup with $0$ and with nonzero idempotent $e$. $T$ is a rectangular $0$-subgroup of $S$ containing $e$ if and only if $T=L'(e)H'(e)^0R'(e)$
for some $L'(e) \subset L(e)$ with $e \in L'(e)$, some $R'(e) \subset R(e)$ with $e \in R'(e)$, such that $R'(e)L'(e) \subset \{e\}^0$, $R'(e)x \neq \{0\}$ for all $x \in L'(e)$, and $yL'(e) \neq \{0\}$ for all $y \in R'(e)$. Every nonzero element of $T$ has a unique product decomposition in $L'(e)$, $H'(e)$, and $R'(e)$.
\end{corollary}

\begin{proof}
We only prove the uniqueness of product decomposition because the rest of the corollary follows directly from previous results.

Let $x,s \in L'(e)$, $a,b \in H'(e)$, and $y,t \in R'(e)$ such that $xay=sbt$. Then there exist $w \in L'(e)$ and $z \in R'(e)$ such that $wx \neq 0$ and $yz \neq 0$, which implies $wx=e$ and $yz=e$ because $R'(e)L'(e) \subset \{e\}^0$. Therefore $wxayz=wsbtz$, or $a=wsbtz$.
Then $ws \neq 0$ and $tz \neq 0$, so that $ws=e$ and $tz=e$ because $R'(e)L'(e) \subset \{e\}^0$, which implies $a=b$. Therefore $wxay=wsbt$, which implies $y=t$, and $xayz=sbtz$, which implies $x=s$.
\end{proof}

\end{document}